\newtheorem{theorem}{Theorem}
\newtheorem{lemma}[theorem]{Lemma}
\newtheorem{corollary}[theorem]{Corollary}
\newtheorem{remark}[theorem]{Remark}

\newcommand{\up}{\mathbb{N}^\omega/\mathcal{F}}
\newcommand{\ap}{\mathbb{A}^\omega/\mathcal{F}}
\newcommand{\proof}{\noindent{\bf Proof.}\hspace{1mm}}
\newcommand{\qed}{{$\Box$}\medskip

}
\author{Juliette Kennedy\thanks{Research partially supported by
grant 40734 of the Academy of Finland.}\\
 Department of Mathematics\\
 University of Helsinki\\
 Helsinki, Finland\\
\and 
Saharon Shelah\thanks{Research partially supported by the United States-Israel
        Binational Science Foundation. Publication number [728] }\\
 Institute of Mathematics\\
 Hebrew University\\
 Jerusalem, Israel\\ 
 }
\title{On embedding models of arithmetic of cardinality $\aleph_1$ into
reduced powers}
\documentclass[12pt]{article}
\usepackage{amsfonts}
\usepackage{amssymb}
\usepackage{latexsym}
\begin{document}
\maketitle
\begin{abstract}

In the early 1970's  S.Tennenbaum proved that all countable models of $PA^- + 
\forall_1 -Th(\mathbb{N})$ are embeddable into the reduced product
$\mathbb{N}^\omega/\mathcal{F}$, where $\mathcal{F}$ is the cofinite filter.
In this paper we show that if $M$ is a model of $PA^- + 
\forall_1 -Th(\mathbb{N})$, and $|M|=\aleph_1$, then $M$ is embeddable into
$\mathbb{N}^\omega/D$, where $D$ is any regular filter on $\omega$.

\end{abstract}

\section{Preliminaries}

Let LA, the language of arithmetic, be the first order language with
non-logical symbols
${+,\cdot,0,1,\leq}$. $\mathbb{N}$ denotes the standard LA
structure.
We shall be concerned with the following theories: The theory
$\forall_1$-$Th(\mathbb{N})$,
defined as the set of all universal
formulas true in the standard LA structure $\mathbb{N}$. Henceforth we refer
to theories satisfying $\forall_1$-$Th(\mathbb{N})$ as Diophantine correct.
We will also refer to the theory $PA^-$, which consists of
the following axioms ($x<y$ abbreviates $x\leq y \wedge x\neq y$):

1.) $\forall x,y,z((x+y)+z=x+(y+z))$

2.) $\forall x,y(x+y=y+x)$

3.) $\forall x,y,z((x\cdot y)\cdot z=x\cdot (y\cdot z))$

4.) $\forall x,y(x\cdot y=y\cdot x)$

5.) $\forall x,y,z(x\cdot (y+z)=x\cdot y+x\cdot z)$

6.) $\forall x((x+0=x)\wedge (x\cdot 0=0))$

7.) $\forall x(x\cdot 1=x)$

8.) $\forall x,y,z((x<y \wedge y<z)\rightarrow x<z)$

9.) $\forall x\;x\leq x$

10.) $\forall x,y(x<y \vee x=y \vee y<x)$

11.) $\forall x,y,z(x<y \rightarrow x+z<y+z)$

12.) $\forall x,y,z (0<z\wedge x<y\rightarrow x\cdot z<y\cdot z)$

13.) $\forall x,y(x<y\rightarrow \exists z(x+z=y)$

14.) $0<1 \wedge \forall x(x>0 \rightarrow x\geq 1)$

15.) $\forall x(x\geq0)$\\

Thus the theory $PA^-$ is the theory of nonnegative parts of
discretely ordered rings. For interesting examples of these models
see \cite{k1}.

\section{The Countable Case}

We begin by presenting the two embedding theorems of Stanley Tennenbaum,
which represent countable models of $PA^-$ by means of sequences of real
numbers. We note that Theorem \ref{1} follows from the 
$\aleph_1$-saturation of 
the structure $\mathbb{N}^\omega/\mathcal{F}$ (\cite{jo}), 
however we present Tennenbaum's construction 
as it constructs the embeddings directly:

We consider first the
reduced power (of LA structures) $\mathbb{N}^\omega/\mathcal{F}$, where
$\mathcal{F}$ is the cofinite filter in the boolean algebra of subsets
of $\mathbb{N}$.
Let $\mathbb{A}$ be the standard LA
structure with domain all nonnegative real algebraic numbers.
We also consider the reduced power $\mathbb{A}^\omega/\mathcal{F}$.
  
If $f$ is a function from $\mathbb{N}$ to $\mathbb{N}$, let $[f]$ denote
the
equivalence class of $f$ in $\up$. We use a similar notation for
$\ap$. When no confusion is possible, we will use
$f$ and $[f]$ interchangeably.

\begin{theorem}\label{1}{
 (Tennenbaum)}
Let $M$ be a countable
Diophantine correct model
of $PA^-$. Then $M$ can be embedded in $\up$.
\end{theorem}

\proof
Let $
m_1,m_2,\ldots $ be the distinct elements of $M$.
Let $P_1, P_2, \ldots $ be all
polynomial equations over $\mathbb{N} $ in the variables $x_1,x_2,\ldots
$ such that $M\models P_i(x_1/m_1,x_2/m_2,\ldots )$.
Each system of equations  $P_1\wedge \cdots \wedge P_n$ has
a solution in
$M$. Thus, by Diophantine correctness, there is
a sequence of natural numbers
$v_1(n),v_2(n), \ldots $ for which
\[ \mathbb{N} \models (P_1\wedge \cdots \wedge P_n)
(x_1/v_1(n),x_2/v_2(n),\ldots ).
\]
Note that if the variable $x_i$ does not appear in
$P_1\wedge\cdots\wedge P_n$, then the choice of $v_i(n)$ is
completely arbitrary.
Our embedding $h:M\longrightarrow \up $ is  given by:
\[m_i
\longmapsto
[\lambda n. v_i(n)] . \]

In the figure below, the i-th row is the solution in integers to
$P_1\wedge \cdots \wedge P_n$, and the i-th column ``is''
$h(m_i)$.

\[ \begin{tabular}{l|lllll}
            &$m_1$     &$m_2$    &$\ldots$    &$m_n$     &$\ldots$ \\
\hline
 $P_1$     &$v_1(1)$  &$v_2(1)$ &$\ldots$    &$v_n(1)$  &$\ldots$ \\
 $P_2$      &$v_1(2)$  &$v_2(2)$ &$\ldots$    &$v_n(2)$  &$\ldots$ \\
 $\:\vdots$   &$\;\;\vdots$  &$\;\;\vdots$ &            &$\;\;\vdots$
&\\

 $P_n$      &$v_1(n)$  &$v_2(n)$ &$\ldots$    &$v_n(n)$  &$\ldots$ \\

 $\:\vdots$   &$\;\;\vdots$  &$\;\;\vdots$ &            &$\;\;\vdots$
&\\

\end{tabular} \]

Note that if $m_e$ is the element $0^M$ of $M$, then the polynomial
equation $x_e=0$
appears as one of the $P$'s. It follows that, for $n$ sufficiently
large, $v_e(n) =0$. Thus $h(0^M)$ is the equivalence class of the zero
function. Similarly, $h$ maps every standard integer of $M$ to the class
of the corresponding constant function.

We show that $h$ is a homomorphism. Suppose $M\models
m_i+ m_j = m_k$. Then the polynomial $x_i+x_j=x_k$ must be one of the
$P$'s, say $P_r$. If $n\geq r$, then by construction
$v_i(n)+v_j(n)=v_k(n)$. Hence $\mathbb{N}^\omega/\mathcal{F}\models
h(m_i)+h(m_j)=h(m_k)$, as
required. A similar argument works for multiplication. Suppose $M\models
m_i\leq m_j$. By an axiom of $PA^-$, for some $k$, $M\models
m_i+m_k=m_j$.
Thus, as we have shown, $h(m_i)+h(m_k)=h(m_j)$. It follows from
the definition of the relation $\leq$ in $\mathcal{N}$ that $h(m_i)\leq
 h(m_j)$.

To see that $h$ is one to one, suppose that $m_i\not= m_j$. Since in
models of
$PA^-$ the order relation is total, we may assume that $m_i < m_j$.
Again by the axioms of $PA^-$, we can choose $m_k$ such that
$m_i+m_k+1=m_j$. As we have shown, $\mathbb{N}^\omega/\mathcal{F}\models h(m_i)+h(m_k)+h(1)=h(m_j)$. Since
h(1)
is the class of the constant function 1, it follows that $h(m_i)\neq
h(m_j)$. \qed

\begin{corollary}
Let $M$ be a countable model of the $\forall_1$-$Th(N)$. Then $M$ can be
embedded in $\mathcal{N}$.
\end{corollary}

\proof
The models of the $\forall_1$-$Th(N)$ are precisely the substructures of
models of $Th(N)$.  Thus,
$M$ extends to a
model of $PA^-$, which can be embedded in $\mathcal{N}$ as in Theorem \ref{1}.
\qed

Before turning to the theorem for the non-Diophantine correct
case, we
observe first that the given embedding depends upon a particular choice
of enumeration $m_1,m_2,\ldots $ of $M$, since different enumerations
will
in general produce different polynomials. We also note that different
choices of solution yield different embeddings. Also, as we
shall see below, we need {\it not} restrict ourselves to Diophantine
formulas: we can carry out the construction for LA formulas
of any complexity which hold in $M$.

We state the
non-Diophantine correct case of the theorem:

\begin{theorem}\label{2}{
(Tennenbaum)} Let $M$ be a countable model of
$PA^-$. Then
$M$ can be embedded in $\ap$.
\end{theorem}

\proof Given an enumeration $m_1,m_2,\ldots $ of $M$,
we form conjunctions of polynomial equations $P_n$ exactly
as before. We
wish to produce solutions of $P_1\wedge \cdots \wedge P_n$
in the nonnegative algebraic
reals for each $n$. We proceed as follows: The model $M$
can be embedded in a real closed field $F$ by a standard
construction. (Embed $M$ in an ordered integral domain,
then form the (ordered) quotient field, and then the real closure.)
Choose
$k$ so large that
$x_1,\ldots,x_k$ are all
the  variables that occur in the conjunction $P_1\wedge \cdots \wedge
P_n$. The sentence $\exists x_1\ldots x_k (P_1\wedge \cdots
\wedge P_n \wedge x_1\geq 0\wedge x_2\geq 0\cdots
\wedge x_k \geq 0) $ is true in $M$, hence in F. It is a theorem of
Tarski that the theory
of real closed fields is complete.
 Thus, this same sentence must
be true in the field of real algebraic numbers. This means we can choose
nonnegative algebraic real numbers $v_1(n),v_2(n)\ldots$ satisfying the
conjunction $P_1\wedge \cdots \wedge P_n$. Let
$h:M\longrightarrow
\ap $ be given by \[m_i\longmapsto [\lambda n.v_i(n)]. \]
The
proof that $h$ is a homomorphism, and furthermore an embedding, proceeds
exactly as before, once we note that the equivalence classes all consist
of nonnegative sequences of real algebraic numbers.
\qed

\begin{remark}
Under any of the embeddings given above, if
$M\models PA^-$ then nonstandard elements of $M$ are mapped to
equivalence classes of functions tending to infinity. Why? If $f$ is a
function in the image of $M$, and $f$ does not tend to
infinity, then choose an integer $k$ such that $f$ is less than $k$
infinitely
often. Since $M\models PA^-$, either $[f] \leq [k]  $ or $[k] \leq [f]
$. The second alternative contradicts the definition of $\leq $ in
$\mathcal{N}$. Hence $[f]\leq [k]$, i.e., $[f]$ is standard.
\end{remark}

\begin{remark}
Let $F$ be a countable ordered field. Then
$F$ is embedded in $\mathbb{R} ^{\omega}/\mathcal{F}$, where $\mathbb{R}
$ is the field of real algebraic numbers. The proof is mutatis mutandis
the same as in Theorem \ref{2}, except that due to the presence of negative
elements we must demonstrate differently that the mapping obtained is
one to one. But this must be the case, since every
homomorphism of fields has this property.

\end{remark}

\begin{remark}
For any pair of LA structures $A$ and
$B$ satisfying $PA^-$, if $A$ is countable and if $A$ satisfies the
$\forall_1$-$Th(B)$ then there is an embedding of $A$ into
$B^{\omega}/\mathcal{F}$. In particular, if $M$ is a model of $PA^-$,
then every
countable extension of $M$ satisfying the $\forall_1$-$Th(M) $ can be
embedded in $M^{\omega }/\mathcal{F}$. \end{remark}

\begin{remark}
Given Theorem \ref{1}, one can ask, what is a necessary
  and sufficient condition for a function to belong to a model of arithmetic
  inside  $\mathcal{N}$? For a partial solution to this
  question, see \cite{jk}.
\end{remark}

\section{The Uncountable Case}

We now show that some of the restrictions of Theorem \ref{1} can be to some
extent relaxed, i.e. we will prove Theorem \ref{1} for 
models $M$ of cardinality 
$\aleph_1$ and with an arbitrary regular filter $D$ in 
place of the cofinite filter.
We note that for
 filters \(D\) on
\(\omega\) for which \(B^\omega/D\) is
\(\aleph_1\)-saturated, where \(B\) is the
two element Boolean algebra, this follows
from the result of  Shelah in \(\cite{shb}\),
that the reduced power $\mathbb{N}^\omega/D$  is $\aleph_1$-saturated.

Our strategy is similar to the strategy of the proof of Theorem
\ref{1}, in that we give an inductive proof on larger and larger initial
segments of the elementary diagram of $M$. 
However we must
now consider formulas whose variables are taken from a set of $\aleph_1$
variables. This requires representing each ordinal $\alpha<\omega$ in terms of
finite sets $u^{\alpha}_n$, which sets determine the variables handled at
each stage of the construction. The other technicality we require is the use
of the following function $g(x,y)$, which bounds the size of the formulas
handled at each stage of the induction.

Let $h(n,m)$ = the total number of non-equivalent Diophantine formulas $\phi(x_1, \ldots x_m)$ of
length $\leq n$.
Define 
\begin{eqnarray*}
  g(n,n)&=&h(n,0)\\
g(n,m-1)&=&2+h(g(n,m),m)\cdot (g(n,m)+3)\mbox{, for }m\leq n. \\
\end{eqnarray*}

\begin{theorem}\label{3}
Let $M$ be a model of
$PA^-$ of cardinality \(\aleph_1\) which is Diophantine correct and let $D$ 
be a regular filter on $\omega$. Then
$M$ can be embedded in $\mathbb{N^{\omega}}/D$.
\end{theorem}

\proof
Let \(M=\{a_\alpha:\alpha<\omega_1\}.\) 
 Let $\{A_n\}_{n\in \omega}$ be a family 
witnessing the regularity of $D$.
We define, for each
\(\alpha<\omega_1\)
a function \(f_\alpha\in\mathbb{N}^{\mathbb{N}}\). Our embedding is then
\(a_\alpha\mapsto[ f_\alpha]\). We need first a lemma:

\def\uan{u^{\alpha}_n}
\def\uanp{u^{\alpha}_{n+1}}
\def\ubnp{u^{\beta}_{n+1}}
\def\ubn{u^{\beta}_n}
\def\ugn{u^{\gamma}_n}

\begin{lemma}
There exists a family of sets \(\uan\), with \(\alpha<\omega_1\),
and \(n\in\mathbb{N}\), such that for each \(n,\alpha\)
\begin{description}
\item[(i)] \(|\uan|<n+1\)
\item[(ii)] \(\alpha\in\uan\subseteq\uanp\)
\item[(iii)] \(\bigcup_n\uan=\alpha+1\)
\item[(iv)] \(\beta\in\uan\Rightarrow\ubn=\uan\cap(\beta+1)\)
\item[(v)] \(\lim_{n\rightarrow\infty}\frac{|\uan|}{n+1}=0\)
\end{description}  
\end{lemma}

Suppose we have the lemma and suppose we have defined
\(f_\beta\) for all \(\beta<\alpha\). We choose \(f_\alpha\in
\mathbb{N}^{\mathbb{N}}\)  componentwise, i.e. we choose
\(f_\alpha(n)\) for each \(n\) separately so that \(f_\alpha(n)\)
satisfies the following condition:
\begin{description}
\item[\((*)_{\alpha,n}\):] 
If \(\phi=\phi(...,x_\beta,...)_{\beta\in\uan}\) 
is a Diophantine formula,  such that
the length of \(\phi\) is $\le g(n,|\uan|)$,
then
$$M\models\phi(...,a_\beta,...)_{\beta\in\uan} \Rightarrow
 \mathbb{N}\models\phi(...,f_\beta(n),...)_{\beta\in\uan}.
$$
\end{description}
(A formula $\phi$ is said to be Diophantine if it has the form 
$$\exists x_0,\ldots,\exists
x_{n-1}(t_1(x_0,\ldots,x_{n-1})=t_2(x_0,\ldots,x_{n-1})),$$ where $t_1$ 
and $t_2$ are LA-terms.)
Now suppose 
$\alpha =0, n > 0.$ By (ii) and (iii), 
 $u^0_n =\{0\}$ for all  $n$. We claim that $(*)_{0,n}$
holds for each $n$. To see this, fix $n$ and let

\begin{eqnarray*} 
\Phi &=& \{\phi(x) \mid M\models \phi(a_0),
 \mbox {where $\phi$ is Diophantine with $|\phi| \leq g(n,1)$}\}.\\  
\end{eqnarray*}
This is -- up to equivalence -- a finite set of formulas. Also, $M\models \exists x \bigwedge \{\phi(x) \mid
\phi \in \Phi\}$ and therefore by Diophantine correctness $\mathbb{N} \models  \exists x 
\bigwedge \{\phi(x) \mid \phi \in \Phi\}$. If $k$ witnesses this formula,
set $f_0(n)=k$. 
Clearly now \((*)_{0,n}\) holds.

Now assume $(*)_{\beta,n}$ holds for all $\beta < \alpha$ and for each $n$. We
choose $f_{\alpha}(n)$ for each $n$ as follows. Fix $n <\omega$ and
let 

\begin{eqnarray*}
\Phi&=& \{\phi(x_0,x_1,\ldots x_k) \mid M\models 
\phi(a_\alpha,..., a_\beta,...)_{\beta\in \uan\setminus\{\alpha\}}), \\
&&\mbox{where } \phi \mbox{ is Diophantine
and } |\phi| \leq g(n,|\uan|)\}, \\
\end{eqnarray*}
for $k=|\uan|-1$ (the case that \(|\uan|\le 1\) reduces to the previous
case).
This is again a finite set of formulas, up to equivalence. Now 
$$
M\models\bigwedge \{\phi(a_\alpha,\langle a_\beta
\rangle_{\beta\in\uan\setminus\{\alpha\}})|\phi \in \Phi\}
$$ 
and therefore
$$
M\models \exists x\bigwedge\{ \phi(x,\langle a_\beta
\rangle_{\beta\in\uan\setminus\{\alpha\}})|\phi \in \Phi\}.
$$ 
Let $\phi^{\prime}(x_1,\ldots ,x_{k})\equiv \exists x \bigwedge
\{\phi(x,x_1,\ldots, x_{k})|\phi\in \Phi\}$. Note that

\[|\phi^{\prime}|\le 2+h(g(n,|\uan|),|\uan|)\cdot (g(n,|\uan|) +3)=g(n,|\uan|-1).\]
Now let $\gamma $= max $(\uan \setminus \{\alpha\})$. Then $u^{\gamma}_n =
 \uan\setminus\{\alpha\}$, by $(iv)$ of the lemma. Since \((*)_{\gamma,n}\) 
holds for each $n$ and we know
that $M\models \phi^{\prime}(\langle a_{\beta}
\rangle_{{\beta}\in u^{\gamma}_n})$ then we
 know by the induction hypothesis that 
\(\mathbb{N}\models\phi^{\prime}(\langle f_\beta(n)
\rangle_{\beta\in u^{\gamma}_n})\). 
Thus 
\(\mathbb{N}\models\exists 
x\bigwedge\{\phi(x,...,f_\gamma(n),...)_{\beta\in 
u^{\gamma}_n})|\phi\in\Phi\}\).
Let \(k\) witness this formula and set \(f_\alpha(n)=k\).
Then \((*)_{\alpha,n}\) holds with \(a_\alpha\mapsto
[\langle f_\alpha(n)\rangle]\).

This mapping is an embedding of \(M\) into $\mathbb{N^{\omega}}/D$.
To see this, suppose
$M\models \phi(a_\alpha, a_\beta,a_\gamma)$ where $\phi$ is the formula
\(x_0+x_1=x_2\). Let \(n=n_0\) be large enough
so that \(|\phi| \leq g(n_0,|u^{\delta}_{n_0}|)$, 
where \(\delta\) is chosen so that \(u^{\delta}_{n_0}\)
contains \({\alpha,\beta,\gamma}\). Then for all \(n \geq n_0\), \((*)_
{\delta,n}\) holds and \(\mathbb{N}\models
f_{\alpha}(n)+f_{\beta}(n)=f_{\gamma}
(n)\), for all \(n\geq n_0\). 
Let $A_{n_1}$ be an element of the chosen regular family of
$D$ such that $A_{n_1} \cap \{0,\ldots,n_0\}=\emptyset$. Then since  
\(
f_{\alpha}(n)+f_{\beta}(n)=f_{\gamma}
(n)\), for \(n\geq n_0\), this holds also for $n\in A_{n_1}$. 
The proof that multiplication is preserved is
the same so we omit it.
Finally, as in the countable case, we note that our mapping
is one-to-one and hence an embedding.

We now prove the lemma, by induction on \(\alpha\).
Let $u^0_n=\{0\}$ for all $n<\omega$.
\medskip

\noindent {\bf Case 1. }\(\alpha\) is a successor ordinal, i.e. \(\alpha=\beta+1\).
Let $n_0$ be such that $n\geq n_0$ implies $\frac{|\ubn|}{n} <\frac{1}{2}.$
then we set 

$$
\uan =\left\{\begin{array}{ll}
\{\alpha\} & n<n_0\\
\ubn \cup \{\alpha\} & n\geq n_0.
\end{array} \right.$$ 
Then (i), (ii) and (iii)   are trivial. 
Proof of (iv): Suppose $\gamma \in \uan = \ubn \cup \{\alpha\}.$ 
\medskip

\noindent {\bf 
Case 1.1.} \(\gamma=\alpha\). Then
\begin{eqnarray*}
  \ugn&=&\uan\\
&=&\uan\cap\alpha+1\\
&=&\uan\cap(\gamma+1).
\end{eqnarray*}
\medskip

\noindent {\bf 
Case 1.2. }$\gamma \in \ubn.$ Then by the induction hypothesis
\begin{eqnarray*}
u^\gamma_n &=&\ubn \cap (\gamma +1)\\
&=& \uan \cap (\gamma +1),
\end{eqnarray*}
since $\gamma \in \ubn$ implies $\gamma \leq \beta.$
\medskip

Proof of (v):
\[\lim_{n\rightarrow\infty}\frac{|\uan|}{n}\le
\lim_{n\rightarrow\infty}\frac{|\ubn|+1}{n}=0,\]
by the induction hypothesis.
\medskip

\noindent
{\bf Case 2.} \(\delta\) is a limit ordinal \(>0\). Let
\(\delta_n\) be an increasing cofinal \(\omega\)-sequence
converging to \(\delta\), for all \(\delta<\omega_1\). 
Let us choose natural
numbers \(n_0,n_1,...\) such that \(\delta_i\in u^{\delta_{i+1}}_{n_i}\)
and \(n\ge n_{i+1}\) implies \(n_i\cdot |u_n^{\delta_{i+1}}|
<n\). Now we let
\[\uan=u^{\delta_i}_n\cup\{\alpha\}\mbox{, if }n_i\le n<n_{i+1}.\]
To prove (iv), let \(\gamma\in\uan\). We wish to show that
\(u^\gamma_n=u^\alpha_n\cap(\gamma+1)\). Suppose \(n_i\le n<n_{i+1}\).
Then \(\uan=u^{\delta_i}_n\cup\{\alpha\}.\) But then
\(\gamma\in\uan\) implies \(\gamma\in u^{\delta_i}_n\cup\{\alpha\}\).
If \(\gamma=\alpha\), then as before,
\(u^\gamma_n=u^\alpha_n\cap(\gamma+1)\).
So suppose \(\gamma\in u^{\delta_i}_n\). Then
\[u^{\gamma}_n=u^{\delta_i}\cap(\gamma+1)=
\uan\cap(\gamma+1).\]
Finally, we prove (iii): Let \(\gamma\in u^\delta_n\), \(\delta\) a limit.
Let \(n_i\) be such that \(\uan=u^{\delta_i}_{n_i}\cup\{\alpha\}\).
We have \(u^{\delta_i}_{n_i}=u^{\delta_{i+1}}_{n_i}\cap(\delta_i+1)\).
Therefore \(\gamma\in u^{\delta_{i+1}}_{n_i}\cup\{\delta\}=u^\delta_{n+1}\).
To prove (v), let \(\epsilon>0\) and choose \(i\) so that
\(\epsilon\cdot n_i>2\). We observe that \(n_i\le n\) implies
\[\lim_{n\rightarrow\infty}\frac{|\uan|}{n}=
\lim_{n\rightarrow\infty}\frac{|u^{\delta_i}_n|}{n}+\frac{1}{n}<
\frac{1}{n_i}+\frac{1}{n}<\epsilon.\]



\qed

In \cite{ks} models of higher cardinality are considered, and
embedding theorems are obtained under a set theoretic assumption.

\end{document}